\documentclass[12pt]{article}
\usepackage[psamsfonts]{amssymb}
\usepackage{amsthm,amsmath}

\title{Reliability analysis of semicoherent systems through their lattice polynomial descriptions}

\author{Alexander Dukhovny\\
Mathematics Department, San Francisco State University \\
San Francisco, CA 94132, USA \\
dukhovny[at]math.sfsu.edu %
\and
Jean-Luc Marichal\thanks{Corresponding author.}\\
Mathematics Research Unit, University of Luxembourg \\
162A, avenue de la Fa\"{\i}encerie, L-1511 Luxembourg, Luxembourg \\
jean-luc.marichal[at]uni.lu %
}

\date{Release 1.0, September 8, 2008}

\begin{document}
\maketitle

\newenvironment{disarray}%
 {\everymath{\displaystyle\everymath{}}\array}%
 {\endarray}

\theoremstyle{plain}
\newtheorem{theorem}{Theorem}
\newtheorem{lemma}[theorem]{Lemma}
\newtheorem{proposition}[theorem]{Proposition}
\newtheorem{corollary}[theorem]{Corollary}

\theoremstyle{definition}
\newtheorem{definition}[theorem]{Definition}
\newtheorem{example}[theorem]{Example}

\theoremstyle{remark}
\newtheorem{remark}[theorem]{Remark}

\newcommand{\N}{\mathbb{N}}                     
\newcommand{\Z}{\mathbb{Z}}                     
\newcommand{\R}{\mathbb{R}}                     
\newcommand{\Vspace}{\vspace{2ex}}                  
\newcommand{\Ind}{\mathrm{Ind}}

\begin{abstract}
A semicoherent system can be described by its structure function or, equivalently, by a lattice polynomial function expressing the system
lifetime in terms of the component lifetimes. In this paper we point out the parallelism between the two descriptions and use the natural
connection of lattice polynomial functions and relevant random events to collect exact formulas for the system reliability. We also discuss the
equivalence between calculating the reliability of semicoherent systems and calculating the distribution function of a lattice polynomial
function of random variables.
\end{abstract}


\section{Introduction}

Consider a semicoherent system made up of nonrepairable components. Such a system can be described by its structure function, which expresses at
any time the state of the system in terms of the states of its components. Equivalently, the system can be described by a lattice polynomial
(l.p.)\ function which expresses the system lifetime in terms of the component lifetimes.

In this paper, 
we point out the formal parallelism
between both descriptions, we collect exact formulas for the system reliability, and we show that calculating the reliability of semicoherent
systems is equivalent to calculating the distribution function of an l.p.\ function of random variables.

We also consider the more general case where there are collective upper bounds on lifetimes of certain subsets of components imposed by external
conditions (such as physical properties of the assembly) or even collective lower bounds imposed for instance by back-up blocks with constant
lifetimes. In terms of lifetimes, such systems can be described by weighted lattice polynomial (w.l.p.) functions. In terms of state variables,
we will see that a ``weighted version'' of the structure functions is required.

This paper is organized as follows. In \S{2} we discuss the parallelism between the description of semicoherent systems by structure functions
and by the corresponding l.p.\ functions. In particular, in \S\ref{sec:SystDesc}, Theorem~\ref{thm:SFvsLP} uses the natural connection between
lattice polynomial functions and relevant random events to establish a centrally important relation between the lifetimes of the system and its
components. In \S{3} we yield exact formulas for the system reliability in case of independent arguments and in general. In turn, those formulas
make it possible to provide exact formulation of reliability parameters such as the mean time-to-failure of the system. In \S{4} we generalize
our results by considering lower and upper bounds on lifetimes of certain components. Finally, in \S{5} we examine how our results can supply
exact formulas for the distribution and moments of w.l.p.\ functions of random variables.

For any numbers $\alpha,\beta\in \overline{\R}:=[-\infty,\infty]$ and any subset $A\subseteq [n]:=\{1,\ldots,n\}$, let
$\mathbf{e}^{\alpha,\beta}_A$ denote the characteristic vector of $A$ in $\{\alpha,\beta\}^n$, that is, the $n$-tuple whose $i$th coordinate is
$\beta$, if $i\in A$, and $\alpha$, otherwise. Also, the $L_1$-norm of any binary vector $\mathbf{x}\in\{0,1\}^n$ is denoted
$|\mathbf{x}|=\sum_{i=1}^nx_i$.

\section{Structure function and l.p.\ function}
\label{sec:StrF-LatPol}

In this section we recall the main concepts and results related to structure functions of semicoherent systems. We also point out the
parallelism between the description of a system by its structure function and the description of this system by an l.p.\ function of the
component lifetimes.

\subsection{Structure function}

Consider a system consisting of $n$ components that are interconnected. The \emph{state}\/ of a component $i\in [n]$ can be represented by a
Boolean variable $x_i$ defined as
$$
x_i =
\begin{cases}
1, & \mbox{if component $i$ is functioning,}\\
0, & \mbox{if component $i$ is in a failed state.}
\end{cases}
$$
For simplicity, we also introduce the state vector $\mathbf{x}=(x_1,\ldots,x_n)$.

The \emph{state of the system}\/ is described from the component states through a Boolean function $\phi:\{0,1\}^n\to\{0,1\}$, called the
\emph{structure function}\/ of the system and defined as
$$
\phi(\mathbf{x}) =
\begin{cases}
1, & \mbox{if the system is functioning,}\\
0, & \mbox{if the system is in a failed state.}
\end{cases}
$$

We shall assume throughout that the structure function $\phi$ is nondecreasing (the system is then said to be \emph{semicoherent}) and
nonconstant, this latter condition ensuring that $\phi(\mathbf{0})=0$ and $\phi(\mathbf{1})=1$, where $\mathbf{0}:=(0,\ldots,0)$ and
$\mathbf{1}:=(1,\ldots,1)$. For a background on semicoherent systems, see for instance the monographs by Ramamurthy \cite{Ram90} and Rausand and
H{\o}yland \cite{RauHoy04}.

As a Boolean function, the structure function $\phi$ can also be regarded as a set function $v:2^{[n]}\to\{0,1\}$. The correspondence is
straightforward: We have $v(A)=\phi(\mathbf{e}^{0,1}_A)$ for all $A\subseteq [n]$ and
\begin{equation}\label{eq:PhiForm1}
\phi(\mathbf{x}) = \sum_{A\subseteq [n]}v(A)\prod_{i\in A}x_i\prod_{i\in [n]\setminus A}(1-x_i).
\end{equation}
We shall henceforth make this identification and often write $\phi_v(\mathbf{x})$ instead of $\phi(\mathbf{x})$. Clearly, the structure function
$\phi_v$ is nondecreasing and nonconstant if and only if its underlying set function $v$ is nondecreasing and nonconstant.

We also observe that, being a Boolean function, the function $\phi_v$ has a unique expression as a multilinear polynomial in $n$ variables,
\begin{equation}\label{eq:PhiForm2}
\phi_v(\mathbf{x}) = \sum_{A\subseteq [n]}m_v(A)\prod_{i\in A}x_i
\end{equation}
(see for instance Hammer and Rudeanu~\cite{HamRud68}), where the set function $m_v:2^{[n]}\to\Z$ is the \emph{M\"obius transform}\/ of $v$,
defined by
$$
m_v(A)=\sum_{B\subseteq A}(-1)^{|A|-|B|}\, v(B).
$$

Another concept that we shall often use in this paper is the \emph{dual}\/ of the set function $v$, that is, the set function
$v^*:2^{[n]}\to\{0,1\}$ defined by $v^*(A)=1-v([n]\setminus A)$.

By extending formally the structure function $\phi_v$ to $[0,1]^n$ by linear interpolation, we define the \emph{multilinear extension}\/ of
$\phi_v$ (a concept introduced in game theory by Owen \cite{Owe72}), that is, the multilinear polynomial function $\overline{\phi}_v:[0,1]^n\to
[0,1]$ defined by
\begin{equation}\label{eq:MLE1}
\overline{\phi}_v(\mathbf{x}) = \sum_{A\subseteq [n]}v(A)\prod_{i\in A}x_i\prod_{i\in [n]\setminus A}(1-x_i).
\end{equation}

Now, by combining the concepts of M\"obius transform, dual set function, and even the ``coproduct'' operation $\amalg$, defined by $\amalg_i
x_i=1-\Pi_i(1-x_i)$, we can easily derive various useful forms of the structure function. Each of these forms is a polynomial expression of the
function $\phi_v$ and hence, when formally regarded as a function from $[0,1]^n$ to $[0,1]$, it identifies with the corresponding multilinear
extension $\overline{\phi}_v$; see also Grabisch et al.~\cite{GraMarRou00}. Table~\ref{tab:StrFun} summarizes the best known forms of the
structure function and its multilinear extension.

\begin{table}[tbp]
$$
\begin{disarray}{|c|c|}
\hline \mbox{Name} & \phi_v(\mathbf{x})~\mbox{and}~\overline{\phi}_{v_{\mathstrut}}^{\mathstrut}(\mathbf{x}) \\
\hline 
\mbox{Primal form} & \sum\limits_{A\subseteq [n]}^{\mathstrut} v(A)\prod\limits_{i\in A}x_i\prod\limits_{i\in [n]\setminus A}(1-x_i)\\ &\\
\mbox{Dual form} & 1-\sum\limits_{A\subseteq [n]}v^*(A)\prod\limits_{i\in [n]\setminus A}x_i\prod\limits_{i\in A}(1-x_i)\\ &\\
\mbox{Primal M\"obius form} & \sum\limits_{A\subseteq [n]}m_v(A)\prod\limits_{i\in A}x_i\\ &\\
\mbox{Dual M\"obius form} & \sum\limits_{A\subseteq [n]}m_{v^*}(A)\coprod\limits_{i\in A}x_i\\ &\\
\mbox{Disjunctive normal form} & \coprod\limits_{A\subseteq [n]} v(A)\prod\limits_{i\in A} x_i\\ &\\
\mbox{Conjunctive normal form} & \prod\limits_{A\subseteq [n]_{\mathstrut}} v^*(A)\coprod\limits_{i\in A} x_i\\
\hline
\end{disarray}
$$
\caption{Various forms of the structure function and its multilinear extension} \label{tab:StrFun}
\end{table}

\subsection{L.p.\ function}

For any event $E$, let $\mathrm{Ind}(E)$ represent the \emph{indicator random variable}\/ that gives $1$ if $E$ occurs and $0$ otherwise. For
any $i\in [n]$, we denote by $T_i$ the random \emph{time-to-failure}\/ of component $i$ and we denote by $X_i(t)=\Ind(T_i>t)$ the random
\emph{state at time $t\geqslant 0$}\/ of component $i$. For simplicity, we introduce the random time-to-failure vector
$\mathbf{T}=(T_1,\ldots,T_n)$ and the random state vector $\mathbf{X}(t)=(X_1(t),\ldots,X_n(t))$ at time $t\geqslant 0$. We also denote by $T_S$
the random time-to-failure of the system and by $X_S(t)=\Ind(T_S>t)$ the random state at time $t\geqslant 0$ of the system.

The structure function $\phi_v$ clearly induces a functional relationship between the variables $T_1,\ldots,T_n$ and the variable $T_S$. As we
will see in Theorem~\ref{thm:SFvsLP}, $T_S$ is always an l.p.\ function of the variables $T_1,\ldots,T_n$. Just as for the structure function,
this l.p.\ function provides a complete description of the structure of the system.

Let us first recall the concept of l.p.\ function of real variables; see for instance Birkhoff~\cite[\S{II.5}]{Bir67} and
Gr\"atzer~\cite[\S{I.4}]{Grae03}. Let $L\subseteq \overline{\R}$ denote a totally ordered bounded lattice whose lattice operations $\wedge$ and
$\vee$ are respectively the minimum and maximum operations. Denote also by $a$ and $b$ the bottom and top elements of $L$.

\begin{definition}
The class of \emph{lattice polynomial}\/ (l.p.) functions from $L^n$ to $L$ is defined as follows:
\begin{enumerate}
\item[(i)] For any $k\in [n]$, the projection $\mathbf{t}\mapsto t_k$ is an l.p.\ function from $L^n$ to $L$.

\item[(ii)] If $p$ and $q$ are l.p.\ functions from $L^n$ to $L$, then $p\wedge q$ and $p\vee q$ are l.p.\ functions from $L^n$ to $L$.

\item[(iii)] Every l.p.\ function from $L^n$ to $L$ is constructed by finitely many applications of the rules (i) and (ii).
\end{enumerate}
\end{definition}

Clearly, any l.p.\ function $p:L^n\to L$ is nondecreasing and nonconstant. Furthermore, it was proved (see for instance
Birkhoff~\cite[\S{II.5}]{Bir67}) that such a function can be expressed in \emph{disjunctive}\/ and \emph{conjunctive}\/ normal forms, that is,
there always exist nonconstant set functions $w^d:2^{[n]}\to\{a,b\}$ and $w^c:2^{[n]}\to\{a,b\}$, with $w^d(\varnothing)=a$ and
$w^c(\varnothing)=b$, such that
\begin{equation}\label{eq:LPDisjConj}
p(\mathbf{t})=\bigvee_{\textstyle{A\subseteq [n]\atop w^d(A)=b}}\bigwedge_{i\in A} t_i = \bigwedge_{\textstyle{A\subseteq [n]\atop
w^c(A)=a}}\bigvee_{i\in A} t_i.
\end{equation}
Clearly, the set functions $w^d$ and $w^c$ that disjunctively and conjunctively define the polynomial function $p(\mathbf{t})$ in
(\ref{eq:LPDisjConj}) are not unique. However, it can be shown \cite{Marc} that, from among all the possible set functions that disjunctively
define $p(\mathbf{t})$, only one is nondecreasing. Similarly, from among all the possible set functions that conjunctively define
$p(\mathbf{t})$, only one is nonincreasing. These special set functions are given by
$$
w^d(A)=p(\mathbf{e}^{a,b}_A)\quad\mbox{and}\quad w^c(A)=p(\mathbf{e}^{a,b}_{[n]\setminus A}).
$$
The l.p.\ function disjunctively defined by a given nondecreasing set function $w:2^{[n]}\to\{a,b\}$ will henceforth be denoted $p_w$. We then
have
$$
p_w(\mathbf{t})=\bigvee_{\textstyle{A\subseteq [n]\atop w(A)=b}}\bigwedge_{i\in A} t_i = \bigwedge_{\textstyle{A\subseteq [n]\atop
w^*(A)=b}}\bigvee_{i\in A} t_i,
$$
where $w^*$ is the \emph{dual}\/ of $w$, defined as
\begin{equation}\label{eq:wDual}
w^*=\gamma\circ v^*,
\end{equation}
the function $\gamma:\{0,1\}\to\{a,b\}$ being a simple transformation defined by $\gamma(0)=a$ and $\gamma(1)=b$.

\subsection{System descriptions}
\label{sec:SystDesc}

The following theorem points out the one-to-one correspondence between the structure function and the l.p.\ function that expresses $T_S$ in
terms of the variables $T_1,\ldots,T_n$. As lifetimes are $[0,\infty]$-valued, we shall henceforth assume without loss of generality that
$L=[0,\infty]$, that is, $a=0$ and $b=\infty$. We also make use of the transformation $\gamma$ as defined in (\ref{eq:wDual}).

\begin{theorem}\label{thm:SFvsLP}
Consider a system whose structure function $\phi_v:\{0,1\}^n\to\{0,1\}$ is nondecreasing and nonconstant. Then we have
\begin{equation}\label{eq:TpT}
T_S=p_w(T_1,\ldots,T_n),
\end{equation}
where $w=\gamma\circ v$. Conversely, any system fulfilling (\ref{eq:TpT}) for some l.p.\ function $p_w:L^n\to L$ has the nondecreasing and
nonconstant structure function $\phi_v$, where $v=\gamma^{-1}\circ w$.
\end{theorem}

\begin{proof}
The proof mainly lies on the distributive property of the indicator function $\Ind(\cdot)$ with respect to disjunction and conjunction, namely
$$
\Ind(E\vee E') = \Ind(E)\vee\Ind(E') \quad\mbox{and}\quad\Ind(E\wedge E') = \Ind(E)\wedge\Ind(E')
$$
for any events $E$ and $E'$. Thus, for any $t\geqslant 0$ we have
\begin{eqnarray*}
\Ind(p_w(\mathbf{T})>t) &=& \Ind\bigg(\bigvee_{\textstyle{A\subseteq [n]\atop v(A)=1}}\bigwedge_{i\in A} T_i>t\bigg)\\
&=& \bigvee_{\textstyle{A\subseteq [n]\atop v(A)=1}}\bigwedge_{i\in A} \Ind(T_i>t) \: =\: \coprod_{\textstyle{A\subseteq [n]\atop v(A)=1}}\prod_{i\in A} X_i(t)\\
&=& \phi_v(\mathbf{X}(t)).
\end{eqnarray*}
Hence, we have $T_S=p_w(\mathbf{T})$ if and only if $X_S(t)=\phi_v(\mathbf{X}(t))$ for all $t\geqslant 0$, which completes the proof.
\end{proof}

\begin{remark}
Since $\phi_v$ is a Boolean function, we can always replace in its expression each product $\Pi$ and coproduct $\amalg$ with the minimum
$\wedge$ and the maximum $\vee$, respectively. Thus, Theorem~\ref{thm:SFvsLP} essentially states that $\phi_v$ is also an l.p.\ function that
has just the same max-min form as $p_w$ but applied to binary arguments. More precisely, $\phi_v$ is \emph{similar} to $p_w$ in the sense that
$\gamma\circ\phi_v=p_w\circ(\gamma,\ldots,\gamma)$.
\end{remark}

We observe that many properties of the structure functions can be derived straightforwardly from the properties of the corresponding l.p.\
functions. Let us examine some of them (see for instance Rausand and H{\o}yland \cite[\S{3.11}]{RauHoy04}):

\begin{enumerate}
\item \textbf{Boundary conditions.} From the idempotency of $p$ (that is, $p(t,\ldots,t)=t$ for all $t\in L$), we immediately retrieve the
idempotency of $\phi$, that is, the boundary conditions $\phi(\mathbf{0})=0$ and $\phi(\mathbf{1})=1$.

\item \textbf{Internality.} The internality property of $p$, namely
$$
\bigwedge_{i=1}^n t_i\leqslant p(\mathbf{t})\leqslant \bigvee_{i=1}^n t_i,
$$
corresponds to the following internality property of $\phi$:
$$
\prod_{i=1}^n x_i\leqslant \phi(\mathbf{x})\leqslant \coprod_{i=1}^n x_i.
$$
Note that, in both cases, internality results immediately from increasing monotonicity and idempotency. For instance, we have
$$
\bigwedge_{i=1}^n t_i=p\Big(\bigwedge_{i=1}^n t_i,\ldots,\bigwedge_{i=1}^n t_i\Big)\leqslant p(\mathbf{t})\leqslant p\Big(\bigvee_{i=1}^n
t_i,\ldots,\bigvee_{i=1}^n t_i\Big)=\bigvee_{i=1}^n t_i.
$$

\item \textbf{Pivotal decomposition.} Consider the following median-based decomposition formula \cite{Marc}, which holds for any l.p.\ function:
\begin{equation}\label{eq:MbDF}
p(\mathbf{t})=\mathrm{median}\big(p(a_i,\mathbf{t}),\, t_i,\, p(b_i,\mathbf{t})\big),
\end{equation}
where the ternary median function is defined as
$$
\mathrm{median}(x_1,x_2,x_3):=(x_1\wedge x_2)\vee (x_2\wedge x_3)\vee (x_3\wedge x_1),
$$
and where $(a_i,\mathbf{t})$ (resp.\ $(b_i,\mathbf{t})$) represents the vector $\mathbf{t}$ whose $i$th coordinate has been replaced with $a$
(resp.\ $b$). From this formula we derive the following property of the structure function:
\begin{eqnarray*}
\phi(\mathbf{x}) &=& \mathrm{median}\big(\phi(0_i,\mathbf{x}),\, x_i,\, \phi(1_i,\mathbf{x})\big)\\
&=& \phi(0_i,\mathbf{x}) \vee \big(x_i \wedge \phi(1_i,\mathbf{x})\big)\: =\: \phi(0_i,\mathbf{x}) \amalg \big(x_i \,\Pi\,\phi(1_i,\mathbf{x})\big)\\
&=& \phi(0_i,\mathbf{x})+ x_i\phi(1_i,\mathbf{x})-x_i\phi(0_i,\mathbf{x})\phi(1_i,\mathbf{x})
\end{eqnarray*}
and hence we retrieve the \emph{pivotal decomposition}\/ of the structure function, namely
$$
\phi(\mathbf{x})=x_i\,\phi(1_i,\mathbf{x})+ (1-x_i)\phi(0_i,\mathbf{x}).
$$

\item \textbf{Structures represented by paths and cuts.} From any nonconstant and nondecreasing set function $w:2^{[n]}\to\{a,b\}$, define the
set function $u_w:2^{[n]}\to\{a,b\}$ as
\begin{equation}\label{eq:MinReprUw}
u_w(A)=
\begin{cases}
b, & \mbox{if $w(A)=b$ and $w(B)=a$ for all $B\varsubsetneq A$,}\\
a, & \mbox{otherwise}.
\end{cases}
\end{equation}
The disjunctive and conjunctive representations of the l.p.\ function $p_w$ having a minimal number of terms write (see Marichal
\cite[Proposition~8]{Marc})
\begin{equation}\label{eq:LPMinRep}
p_w(\mathbf{t})=\bigvee_{\textstyle{A\subseteq [n]\atop u_w(A)=b}}\bigwedge_{i\in A} t_i = \bigwedge_{\textstyle{A\subseteq [n]\atop
u_{w^*}(A)=b}}\bigvee_{i\in A} t_i.
\end{equation}
Let us show that these representations are in one-to-one correspondence with the representations of the structure function by minimal paths and
cuts. Recall that a \emph{path set}\/ $P\subseteq [n]$ is a set of components which by functioning ensures that the system is functioning.
Similarly, a \emph{cut set}\/ $K\subseteq [n]$ is a set of components which by failing causes the system to fail. In other terms, $P\subseteq
[n]$ is a path set if $v(P)=1$ and $K\subseteq [n]$ is a cut set if $v([n]\setminus K)=0$. A path (resp.\ cut) set is \emph{minimal}\/ if it
does not contain any proper path (resp.\ cut) set.

It is known that if $P_1,\ldots,P_r$ are the minimal path sets and $K_1,\ldots,K_s$ are the minimal cut sets, then
$$
\phi_v(\mathbf{x})=\coprod_{j=1}^r\prod_{i\in P_j}x_i = \prod_{j=1}^s\coprod_{i\in K_j}x_i.
$$
The corresponding formulas for the l.p.\ function write
$$
p_w(\mathbf{t})=\bigvee_{j=1}^r\bigwedge_{i\in P_j}t_i = \bigwedge_{j=1}^s\bigvee_{i\in K_j}t_i
$$
and are exactly the ``minimal'' representations (\ref{eq:LPMinRep}) of $p_w$.

\item \textbf{Extra component connected in series or parallel.} Any l.p.\ function $p:L^n\to L$ fulfills trivially the following functional
equations
\begin{eqnarray*}
p(u\wedge t_1,\ldots,u\wedge t_n) &=& u\wedge p(t_1,\ldots,t_n), \\
p(u\vee t_1,\ldots,u\vee t_n) &=& u\vee p(t_1,\ldots,t_n),
\end{eqnarray*}
for arbitrary $u\in L$. These equations mean that connecting in series (resp.\ in parallel) any extra component to the system amounts to
connecting that component in series (resp.\ in parallel) to each component of the system. The corresponding equations for the structure function
are clear. We have
\begin{eqnarray*}
\phi(y\, x_1,\ldots,y\, x_n) &=& y\,\phi(x_1,\ldots,x_n), \\
\phi(y\amalg x_1,\ldots,y\amalg x_n) &=& y\amalg\phi(x_1,\ldots,x_n),
\end{eqnarray*}
for arbitrary $y\in\{0,1\}$.

\item \textbf{Dual structure.} Recall that the dual structure function of a structure function $\phi_v$ is defined as
$\phi_v^D(\mathbf{x})=1-\phi_v(\mathbf{1}-\mathbf{x})$. From this definition, we derive immediately $\phi^D_v=\phi_{v^*}$, and hence from
(\ref{eq:PhiForm1}) we immediately retrieve the dual form of $\phi_v$ (i.e., the second expression in Table~\ref{tab:StrFun}). Using the dual
set function $w^*$ of $w$, as defined in (\ref{eq:wDual}), we see that the corresponding l.p.\ function is the dual of $p_w$, namely
$p^D_w=p_{w^*}$.
\end{enumerate}

\section{Exact reliability calculation}

The \emph{reliability function}\/ of component $i$ is defined, for any $t\geqslant 0$, by
$$
R_i(t)=\Pr(T_i>t)=\Pr(X_i(t)=1)=\mathrm{E}[X_i(t)],
$$
that is, the probability that component $i$ does not fail in the time interval $[0,t]$. Similarly, for any $t\geqslant 0$, the system
reliability function is
$$
R_S(t)=\Pr(T_S>t)=\Pr(X_S(t)=1)=\mathrm{E}[X_S(t)],
$$
that is, the probability that the system does not fail in the time interval $[0,t]$.

The \emph{mean time-to-failure of component $i$}\/ is defined as $\mathrm{MTTF}_i=\mathrm{E}[T_i]$ and similarly the \emph{mean time-to-failure
of the system}\/ is defined as $\mathrm{MTTF}_S=\mathrm{E}[T_S]$. These expected values can be calculated by the following formulas (see for
instance Rausand and H{\o}yland \cite[\S{2.6}]{RauHoy04})
$$
\mathrm{MTTF}_i=\int_a^b R_i(t)\, dt \quad\mbox{and}\quad \mathrm{MTTF}_S=\int_a^b R_S(t)\, dt.
$$

In this section we yield the main known formulas for the system reliability function in the general case of dependent failures and in the
special case of independent failures. We also provide some additional useful formulas.

\subsection{Dependent failures}
\label{sec:DepF-LP}

Dukhovny~\cite{Duk07} found simple and concise formulas for the system reliability function in case of generally dependent variables
$T_1,\ldots,T_n$. We present them in the following theorem and we provide a shorter proof.


\begin{theorem}\label{thm:RSv}
We have
\begin{eqnarray}
R_S(t) &=& \sum_{A\subseteq [n]} v(A)\,\Pr(\mathbf{X}(t)=\mathbf{e}^{0,1}_A),\label{eq:RSv}\\
R_S(t) &=& 1-\sum_{A\subseteq [n]} v^*(A)\,\Pr(\mathbf{X}(t)=\mathbf{e}^{0,1}_{[n]\setminus A})\label{eq:RSvs}.
\end{eqnarray}
\end{theorem}

\begin{proof}
By (\ref{eq:PhiForm1}), we have
\begin{eqnarray}
R_S(t) = \mathrm{E}[\phi_v(\mathbf{X}(t))] &=& \sum_{A\subseteq [n]}v(A)\,\mathrm{E}\Big[\prod_{i\in A}X_i(t)\prod_{i\in
[n]\setminus A}(1-X_i(t))\Big]\label{eq:RSv1}\\
&=& \sum_{A\subseteq [n]} v(A)\,\Pr(\mathbf{X}(t)=\mathbf{e}^{0,1}_A),\nonumber
\end{eqnarray}
which proves (\ref{eq:RSv}). Formula (\ref{eq:RSvs}) can be proved similarly by using the dual form of $\phi_v$ (i.e., the second expression in
Table~\ref{tab:StrFun}).
\end{proof}

Consider the \emph{joint distribution function}\/ and the \emph{joint survival function}, defined respectively as
$$
F(\mathbf{t}) = \Pr(T_i\leqslant t_i\;\forall i\in [n])\quad\mbox{and}\quad R(\mathbf{t}) = \Pr(T_i>t_i\;\forall i\in [n]).
$$
By using the same argument as in the proof of Theorem~\ref{thm:RSv}, we obtain two further equivalent expressions of $R_S(t)$.

\begin{theorem}\label{thm:RSmvs}
We have
\begin{eqnarray}
R_S(t) &=& \sum_{A\subseteq [n]} m_v(A) \, R(\mathbf{e}^{a,t}_A)\label{eq:RSmv}\\
R_S(t) &=& \sum_{A\subseteq [n]} m_{v^*}(A) \, \big(1-F(\mathbf{e}^{t,b}_{[n]\setminus A})\big)\: =\: 1-\sum_{A\subseteq [n]} m_{v^*}(A) \,
F(\mathbf{e}^{t,b}_{[n]\setminus A}).\nonumber
\end{eqnarray}
\end{theorem}

\begin{proof}
By (\ref{eq:PhiForm2}), we have
\begin{eqnarray*}
R_S(t) = \mathrm{E}[\phi_v(\mathbf{X}(t))] &=& \sum_{A\subseteq [n]} m_v(A) \,\mathrm{E}\Big[\prod_{i\in A}X_i(t)\Big]\\
&=& \sum_{A\subseteq [n]} m_v(A) \,\Pr(T_i>t\;\forall i\in A)\\
&=& \sum_{A\subseteq [n]} m_v(A) \,R(\mathbf{e}^{a,t}_A).
\end{eqnarray*}
Similarly, using the dual M\"obius form of $\phi_v$ (i.e., the fourth expression in Table~\ref{tab:StrFun}), we have
\begin{eqnarray*}
R_S(t) = \mathrm{E}[\phi_v(\mathbf{X}(t))] &=& \sum_{A\subseteq [n]} m_{v^*}(A) \,\mathrm{E}\Big[\coprod_{i\in A}X_i(t)\Big]\\
&=& \sum_{A\subseteq [n]} m_{v^*}(A) \,\big(1-\Pr(T_i\leqslant t\;\forall i\in A)\big)\\
&=& \sum_{A\subseteq [n]} m_{v^*}(A) \,\big(1-F(\mathbf{e}^{t,b}_{[n]\setminus A})\big),
\end{eqnarray*}
and for the last formula, we use the fact that $\sum_{A\subseteq [n]} m_{v^*}(A) = \phi_{v^*}(\mathbf{1})=1$.
\end{proof}

It is noteworthy that Theorem~\ref{thm:RSmvs} immediately provides concise expressions for the mean time-to-failure of the system, namely
\begin{eqnarray*}
\mathrm{MTTF}_S &=& \sum_{A\subseteq [n]} m_v(A) \, \int_0^{\infty} R(\mathbf{e}^{0,t}_A)\, dt,\\
\mathrm{MTTF}_S &=& \sum_{A\subseteq [n]} m_{v^*}(A) \, \int_0^{\infty} \big(1-F(\mathbf{e}^{t,\infty}_{[n]\setminus A})\big)\, dt.
\end{eqnarray*}

Theorem~\ref{thm:RSmvs} may suggest that the complete knowledge of the joint survival (or joint distribution) function is needed for the
calculation of the system reliability function. Actually, as Theorem~\ref{thm:RSv} shows, all the needed information is encoded in the
distribution of the indicator vector $\mathbf{X}(t)$. In turn, the distribution of $\mathbf{X}(t)$ can be easily expressed (see
Dukhovny~\cite{Duk07} and Dukhovny and Marichal \cite{DukMar08}) in terms of the \emph{joint probability generating function} of
$\mathbf{X}(t)$, which is defined by
$$
G(\mathbf{z},t) = \mathrm{E}\Big[\prod_{i=1}^n z_i^{X_i(t)}\Big] \qquad (|z_i|\leqslant 1, \, i\in [n]).
$$
As it is well known, the joint probability generating function has the advantage of being an expectation and yields not only the probabilities
alone but also all kinds of moments via derivatives.

By definition, we have
\begin{equation}\label{eq:jpgf}
G(\mathbf{z},t) = \sum_{\mathbf{x}\in\{0,1\}^n}\Pr(\mathbf{X}(t)=\mathbf{x})\prod_{i\in [n]}z_i^{x_i} = \sum_{A\subseteq [n]}
\Pr(\mathbf{X}(t)=\mathbf{e}_A^{0,1})\prod_{i\in A}z_i
\end{equation}
and hence $G(\mathbf{z},t)$ is a multilinear polynomial in $z_1,\ldots,z_n$, which can be rewritten as
$$
G(\mathbf{z},t) = \sum_{A\subseteq [n]} G(\mathbf{e}_A^{0,1},t) \prod_{i\in A}z_i\prod_{i\in [n]\setminus A}(1-z_i).
$$
Moreover, we can easily show \cite{Duk07,DukMar08} that $G(\mathbf{e}_A^{0,1},t)=F(\mathbf{e}_A^{t,b})$.

On the other hand, from (\ref{eq:jpgf}) it follows that
$$
G(\mathbf{e}_A^{0,1},t) = \sum_{B\subseteq A} \Pr(\mathbf{X}(t)=\mathbf{e}_B^{0,1})
$$
which shows that the set function $A\mapsto \Pr(\mathbf{X}(t)=\mathbf{e}_A^{0,1})$ is the M\"obius transform of the set function $A\mapsto
G(\mathbf{e}_A^{0,1},t)$, that is,
$$
\Pr(\mathbf{X}(t)=\mathbf{e}_A^{0,1}) = \sum_{B\subseteq A} (-1)^{|A|-|B|}\, G(\mathbf{e}_B^{0,1},t).
$$
Combining this latter formula with (\ref{eq:RSv}) enables us to express the system reliability function in terms of $G(\mathbf{z},t)$.

\subsection{Independent failures}

In the case when $T_1,\ldots,T_n$ are independent, which implies that the indicator variables $X_1(t),\ldots,X_n(t)$ are independent for all
$t\geqslant 0$, from (\ref{eq:RSv1}) we obtain the well-known formula
\begin{equation}\label{eq:RsRiInd}
R_S(t) = \sum_{A\subseteq [n]} v(A) \prod_{i\in A}R_i(t)\prod_{i\in [n]\setminus A}(1-R_i(t)).
\end{equation}
Combining (\ref{eq:MLE1}) and (\ref{eq:RsRiInd}), we immediately retrieve the following classical formula (see for instance Rausand and
H{\o}yland \cite[\S{4.5}]{RauHoy04})
$$
R_S(t) = \overline{\phi}_v(R_1(t),\ldots,R_n(t))
$$
and so both $R_S(t)$ and $\mathrm{MTTF}_S$ can be expressed in different forms, according to the expressions of $\overline{\phi}_v$ chosen in
Table~\ref{tab:StrFun}. For instance, using the primal M\"obius form of $\overline{\phi}_v$, we obtain
\begin{eqnarray}
R_S(t) &=& \sum_{A\subseteq [n]} m_v(A) \, \prod_{i\in A} R_i(t),\nonumber\\
\mathrm{MTTF}_S &=& \sum_{A\subseteq [n]} m_v(A) \, \int_0^{\infty} \prod_{i\in A} R_i(t)\, dt.\label{eq:MTTFsMobInd}
\end{eqnarray}


\subsection{Some examples}
\label{sec:ex5}

Let us now examine some typical examples by considering both their structure functions and the corresponding l.p.\ functions:
\begin{enumerate}
\item \emph{Series structure.} If all the components are wired in series, we have $\phi_v(\mathbf{x})=\prod_i x_i$ and
$p_w(\mathbf{t})=\bigwedge_i t_i$. In this case, $w(A)=b$ if and only if $A=[n]$. More generally, we can show \cite[\S{5.4}]{Mar00c} that any
l.p.\ function $p_w$ fulfilling the functional equation
\begin{equation}\label{eq:Minitivity}
p_w(t_1\wedge t'_1,\ldots,t_n\wedge t'_n)=p_w(t_1,\ldots,t_n)\wedge p_w(t'_1,\ldots,t'_n)
\end{equation}
is of the form $p_w(\mathbf{t})=\bigwedge_{i\in B} t_i$ for some subset $B\subseteq [n]$. It then corresponds to a serially connected segment of
components.

The reliability of a series structure with $n$ elements is given by $R_S(t)=\Pr(\wedge_{i=1}^n T_i>t)$. Using Theorems~\ref{thm:RSv} and
\ref{thm:RSmvs}, we also have
$$
R_S(t) = \Pr(\mathbf{X}(t)=\mathbf{1}) = R(\mathbf{e}^{0,t}_{[n]}) = \Pr(T_1>t,\ldots,T_n>t).
$$

\item \emph{Parallel structure.} If all the components are wired in parallel, we have $\phi_v(\mathbf{x})=\coprod_i x_i$ and
$p_w(\mathbf{t})=\bigvee_i t_i$. In this case, $w(A)=b$ if and only if $A\neq\varnothing$. Similarly to the series structures, we can show that
any l.p.\ function $p_w$ fulfilling the functional equation
\begin{equation}\label{eq:Maxitivity}
p_w(t_1\vee t'_1,\ldots,t_n\vee t'_n)=p_w(t_1,\ldots,t_n)\vee p_w(t'_1,\ldots,t'_n)
\end{equation}
is of the form $p_w(\mathbf{t})=\bigvee_{i\in B} t_i$ for some subset $B\subseteq [n]$. It then corresponds to a subsystem of parallel
components.

The reliability of a parallel structure with $n$ elements is given by $R_S(t)=\Pr(\vee_{i=1}^n T_i>t)$. Using Theorems~\ref{thm:RSv} and
\ref{thm:RSmvs}, we also have
$$
R_S(t) = 1-\Pr(\mathbf{X}(t)=\mathbf{0}) = 1-F(\mathbf{e}^{t,\infty}_{\varnothing}) = 1-\Pr(T_1\leqslant t,\ldots,T_n\leqslant t).
$$

\item \emph{$k$-out-of-$n$ structure, for some $k\in [n]$.} By definition, a $k$-out-of-$n$ structure is characterized by the structure function
$$
\phi_v(\mathbf{x}) =
\begin{cases}
1, & \mbox{if $\sum_i x_i\geqslant k$,}\\
0, & \mbox{if $\sum_i x_i< k$.}
\end{cases}
$$
It is then easy to show that
$$
\phi_v(\mathbf{x}) = \coprod_{\textstyle{A\subseteq [n]\atop |A|=k}}\prod_{i\in A}x_i ~ = ~ \prod_{\textstyle{A\subseteq [n]\atop
|A|=n-k+1}}\coprod_{i\in A}x_i,
$$
which, in turn, entails
$$
p_w(\mathbf{t}) = \bigvee_{\textstyle{A\subseteq [n]\atop |A|=k}}\bigwedge_{i\in A}t_i ~ = ~ \bigwedge_{\textstyle{A\subseteq [n]\atop
|A|=n-k+1}}\bigvee_{i\in A}t_i ~ = ~ f_{n-k+1}(\mathbf{t}),
$$
where, for any $k\in [n]$, $f_k:L^n\to L$ is the $k$th order statistic function (see for instance Ovchinnikov \cite{Ovc96}). We recall
\cite[\S{5.5}]{Mar00c} that the $n$ order statistic functions are exactly those l.p.\ functions that are symmetric in their variables. It
follows immediately that a structure is of $k$-out-of-$n$ type for some $k\in [n]$ if and only if its system lifetime is a symmetric function
(which is $f_{n-k+1}$) of the component lifetimes. In this case, $w(A)=b$ if and only if $|A|\geqslant k$, which means that the system is
functioning if at least $k$ components are functioning. Clearly, the minimal representation (\ref{eq:MinReprUw}) is such that $u_w(A)=b$ if and
only $|A|=k$. We also observe that
\begin{eqnarray}
m_v(A) &=&
\begin{cases}
(-1)^{|A|-k}{|A|-1\choose k-1}, & \mbox{if $|A|\geqslant k$},\\
0, & \mbox{otherwise,}
\end{cases}\label{eq:Mobnoon}\\
m_{v^*}(A) &=&
\begin{cases}
(-1)^{|A|-n+k-1}{|A|-1\choose n-k}, & \mbox{if $|A|\geqslant n-k+1$},\\
0, & \mbox{otherwise.}
\end{cases}\nonumber
\end{eqnarray}
The reliability of a $k$-out-of-$n$ structure is given by $R_S(t)=\Pr(f_{n-k+1}(\mathbf{T})>t)$. Using Theorem~\ref{thm:RSv}, we also have
\begin{equation}\label{eq:SymPP}
R_S(t) = \sum_{j=k}^n \Pr(|\mathbf{X}(t)|=j) = \Pr(|\mathbf{X}(t)|\geqslant k),
\end{equation}
and we can show \cite{Duk07,DukMar08} that $\Pr(|\mathbf{X}(t)|=j)=[x^j]G(x\mathbf{1},t)$ is the coefficient of $x^j$ in the $n$th degree
polynomial $G(x\mathbf{1},t)$. On the other hand, combining (\ref{eq:RSmv}) and (\ref{eq:Mobnoon}) gives
$$
R_S(t) =\sum_{j=k}^n (-1)^{j-k}{j-1\choose k-1}\,\sum_{\textstyle{A\subseteq [n]\atop |A|=j}}R(\mathbf{e}^{0,t}_A).
$$
\end{enumerate}

\begin{example}
When $R_i(t)=e^{-\lambda_i t}$ $(i=1,\ldots,n)$, it is convenient to calculate $\mathrm{MTTF}_S$ by using formula (\ref{eq:MTTFsMobInd}).
Indeed, in that case, setting $\lambda_A=\sum_{i\in A}\lambda_i$, we simply obtain (see \cite{Mar06})
\begin{equation}\label{eq:MTTFSExp}
\mathrm{MTTF}_S = \sum_{A\subseteq [n]}m_v(A) \int_0^{\infty} e^{-\lambda_A t}\, dt = \sum_{\textstyle{A\subseteq [n]\atop
A\neq\varnothing}}\frac{m_v(A)}{\lambda_A}.
\end{equation}
Assuming further that the structure is of $k$-out-of-$n$ type, by (\ref{eq:Mobnoon}) we immediately obtain
$$
\mathrm{MTTF}_S = \sum_{\textstyle{A\subseteq [n]\atop |A|\geqslant k}}(-1)^{|A|-k}\,{|A|-1\choose k-1}\,\frac{1}{\lambda_A}\, .
$$
\end{example}

\section{Systems with lower and upper bounds on lifetimes}

Consider now a more general system in which we allow upper and/or lower bounds on lifetimes of certain subsets of components. As shown by
Dukhovny and Marichal \cite{DukMar08}, the structure of such a system can be modelled by means of a w.l.p.\ function, which is an l.p.\ function
constructed from both variables and constants.

\subsection{System descriptions}

We first recall the definition of w.l.p.\ functions (see Goodstein~\cite{Goo67} and Rudeanu~\cite[Chapter~3,\,\S{3}]{Rud01}).

\begin{definition}
The class of \emph{weighted lattice polynomial}\/ (w.l.p.) functions from $L^n$ to $L$ is defined as follows:
\begin{enumerate}
\item[(i)] For any $k\in [n]$ and any $c\in L$, the projection $\mathbf{t}\mapsto t_k$ and the constant function $\mathbf{t}\mapsto c$ are
w.l.p.\ functions from $L^n$ to $L$.

\item[(ii)] If $p$ and $q$ are w.l.p.\ functions from $L^n$ to $L$, then $p\wedge q$ and $p\vee q$ are w.l.p.\ functions from $L^n$ to $L$.

\item[(iii)] Every w.l.p.\ function from $L^n$ to $L$ is constructed by finitely many applications of the rules (i) and (ii).
\end{enumerate}
\end{definition}

It was proved \cite{Goo67} that any w.l.p.\ function $p:L^n\to L$ can be expressed in \emph{disjunctive}\/ and \emph{conjunctive}\/ normal
forms, that is, there exist set functions $w^d:2^{[n]}\to L$ and $w^c:2^{[n]}\to L$ such that
$$
p(\mathbf{t})=\bigvee_{A\subseteq [n]}\Big(w^d(A)\wedge\bigwedge_{i\in A} t_i\Big) = \bigwedge_{A\subseteq [n]}\Big(w^c(A)\vee\bigvee_{i\in A}
t_i\Big).
$$
Moreover, it can be shown \cite{Marc} that, from among all the possible set functions $w^d$ that disjunctively define $p(\mathbf{t})$, only one
is nondecreasing. Similarly, from among all the possible set functions $w^c$ that conjunctively define $p(\mathbf{t})$, only one is
nonincreasing. These special set functions are given by
$$
w^d(A)=p(\mathbf{e}^{a,b}_A)\quad\mbox{and}\quad w^c(A)=p(\mathbf{e}^{a,b}_{[n]\setminus A}).
$$
The w.l.p.\ function defined by a given nondecreasing set function $w:2^{[n]}\to L$ will henceforth be denoted $p_w$.

The following theorem, which generalizes Theorem~\ref{thm:SFvsLP} to w.l.p.\ functions, shows that the system is no longer characterized by a
single structure function but by a one-parameter family of structure functions.

\begin{theorem}\label{thm:SFvsWLP}
With any system fulfilling $T_S=p_w(T_1,\ldots,T_n)$ for some w.l.p.\ function $p_w:L^n\to L$ is associated a unique family of nondecreasing and
nonconstant structure functions $\{\phi_{v_t}: \, t\geqslant 0\}$, where $v_t(A)=\Ind(w(A)>t)$ for all $A\subseteq [n]$, such that
$$
X_S(t)=\phi_{v_t}(\mathbf{X}(t)) \qquad (t\geqslant 0).
$$
\end{theorem}

\begin{proof}
We follow the same reasoning as in the proof of Theorem~\ref{thm:SFvsLP}. For any $t\geqslant 0$, we have
\begin{eqnarray*}
\Ind(p_w(\mathbf{T})>t) &=& \Ind\Big(\bigvee_{A\subseteq [n]}\big(w(A)\wedge\bigwedge_{i\in A} T_i\big)>t\Big)\\
&=& \bigvee_{A\subseteq [n]}\Big(\Ind(w(A)>t)\wedge\bigwedge_{i\in A} \Ind(T_i>t)\Big)\\
&=& \coprod_{A\subseteq [n]}v_t(A)\prod_{i\in A} X_i(t)\\
&=& \phi_{v_t}(\mathbf{X}(t)).
\end{eqnarray*}
Hence, we have $T_S=p_w(\mathbf{T})$ if and only if $X_S(t)=\phi_{v_t}(\mathbf{X}(t))$ for all $t\geqslant 0$, which completes the proof.
\end{proof}

\begin{remark}
According to Theorem~\ref{thm:SFvsWLP}, when modelling systems with collective bounds, it seems much more convenient to use w.l.p.\ functions
rather than families of structure functions. 
\end{remark}

The properties of the family of structure functions can be derived from the properties of the corresponding w.l.p.\ function. Let us examine
some of them:

\begin{enumerate}
\item \textbf{Boundary conditions.} We have $\phi_{v_t}(\mathbf{0})=v_t(\varnothing)=\Ind(w(\varnothing)>t)$ and
$\phi_{v_t}(\mathbf{1})=v_t([n])=\Ind(w([n])>t)$.

\item \textbf{Pivotal decomposition.} The median-based decomposition formula (\ref{eq:MbDF}), which also holds for any w.l.p.\ function, leads
again to the pivotal decomposition of each structure function $\phi_{v_t}$:
$$
\phi_{v_t}(\mathbf{x})=x_i\,\phi_{v_t}(1_i,\mathbf{x})+ (1-x_i)\phi_{v_t}(0_i,\mathbf{x})\qquad (t\geqslant 0).
$$

\item \textbf{Minimal representations.} From a nondecreasing set function $w:2^{[n]}\to L$, define the set functions $u^d_w:2^{[n]}\to L$ and
$u^c_w:2^{[n]}\to L$ as
\begin{eqnarray*}
u^d_w(A)&=&
\begin{cases}
w(A), & \mbox{if $w(B)<w(A)$ for all $B\varsubsetneq A$,}\\
a, & \mbox{otherwise},
\end{cases}
\\
u^c_w(A)&=&
\begin{cases}
w(A), & \mbox{if $w(A)<w(B)$ for all $B\varsupsetneq A$,}\\
b, & \mbox{otherwise}.
\end{cases}
\end{eqnarray*}
The disjunctive and conjunctive representations of the w.l.p.\ function $p_w$ having a minimal number of terms write (see Marichal
\cite[Proposition~8]{Marc})
$$
p_w(\mathbf{t})=\bigvee_{A\subseteq [n]}\Big(u^d_w(A)\wedge\bigwedge_{i\in A} t_i\Big) = \bigwedge_{A\subseteq [n]}\Big(u^c_w([n]\setminus
A)\vee\bigvee_{i\in A} t_i\Big).
$$
The corresponding form of the family of structure functions follows.

\item \textbf{Dual structure.} The dual of a family $\{\phi_{v_t}\, :\, t\geqslant 0\}$ of structure functions is the family $\{\phi_{v^*_t}\,
:\, t\geqslant 0\}$ of structure functions such that
$$
v^*_t(A)=1-\Ind(w([n]\setminus A)>t).
$$
It follows that there is no set function $w^*:2^{[n]}\to L$ such that $v^*_t(A)=\Ind(w^*(A)>t)$.
\end{enumerate}

\begin{example}
Consider the bridge structure as indicated in Figure~\ref{fig:bs} and assume that the time-to-failure of the central component must lie in the
interval $[l,u]$ for some $0\leqslant l\leqslant u\leqslant\infty$. We immediately see that the minimal path sets are $P_1=\{1,4\}$,
$P_2=\{2,5\}$, $P_3=\{1,3,5\}$, and $P_4=\{2,3,4\}$. Hence, the w.l.p.\ function associated with this structure is given by
\begin{eqnarray*}
p_w(\mathbf{t}) &=& (t_1\wedge t_4)\vee (t_2\wedge t_5) \vee (t_1\wedge\mathrm{median}(l,t_3,u)\wedge t_5)\vee (t_2\wedge\mathrm{median}(l,t_3,u)\wedge t_4)\\
&=& (t_1\wedge t_4)\vee (t_2\wedge t_5)\vee (l\wedge t_1\wedge t_5)\vee (u\wedge t_1\wedge t_3\wedge t_5)\vee (l\wedge t_2\wedge t_4)\vee
(u\wedge t_2\wedge t_3\wedge t_4)
\end{eqnarray*}
and the corresponding family of structure functions is
\begin{eqnarray*}
\phi_{v_t}(\mathbf{x}) &=& (x_1x_4)\amalg (x_2x_5)\amalg (\Ind(l>t)x_1x_5)\amalg (\Ind(u>t)x_1x_3x_5)\\
&& \null \amalg (\Ind(l>t)x_2x_4)\amalg (\Ind(u>t)x_2x_3x_4).
\end{eqnarray*}
\end{example}

\setlength{\unitlength}{4ex}
\begin{figure}[htbp]\centering
\begin{picture}(11,4)
\put(3,0.5){\framebox(1,1){$2$}} \put(3,2.5){\framebox(1,1){$1$}} \put(5,1.5){\framebox(1,1){$3$}} \put(7,0.5){\framebox(1,1){$5$}}
\put(7,2.5){\framebox(1,1){$4$}}%
\put(0,2){\line(1,0){1.5}}\put(1.5,2){\line(2,-1){1.5}}\put(5.5,0){\line(-2,1){1.5}}\put(1.5,2){\line(2,1){1.5}}\put(5.5,4){\line(-2,-1){1.5}}%
\put(0,2){\circle*{0.15}}%
\put(9.5,2){\line(1,0){1.5}}\put(5.5,0){\line(2,1){1.5}}\put(9.5,2){\line(-2,-1){1.5}}\put(5.5,4){\line(2,-1){1.5}}\put(9.5,2){\line(-2,1){1.5}}%
\put(11,2){\circle*{0.15}}%
\put(5.5,0){\line(0,1){1.5}}\put(5.5,4){\line(0,-1){1.5}}
\end{picture}
\caption{Bridge structure} \label{fig:bs}
\end{figure}

\subsection{Exact reliability formulas}

Regarding the reliability calculation, Dukhovny and Marichal \cite{DukMar08} established the following result, which is a direct generalization
of Theorem~\ref{thm:RSv}:

\begin{theorem}\label{thm:RSvWLP}
We have
\begin{eqnarray}
R_S(t) &=& \sum_{A\subseteq [n]} v_t(A)\,\Pr(\mathbf{X}(t)=\mathbf{e}^{0,1}_A),\label{eq:RSvWLP}\\
R_S(t) &=& 1-\sum_{A\subseteq [n]} v^*_t(A)\,\Pr(\mathbf{X}(t)=\mathbf{e}^{0,1}_{[n]\setminus A}).\nonumber
\end{eqnarray}
\end{theorem}

Similarly, a direct generalization of Theorem~\ref{thm:RSmvs} is stated in the following theorem:

\begin{theorem}\label{thm:RSmvsWLP}
We have
\begin{eqnarray}
R_S(t) &=& \sum_{A\subseteq [n]} m_{v_t}(A) \, R(\mathbf{e}^{a,t}_A),\label{eq:RSmvWLP}\\
R_S(t) &=& \sum_{A\subseteq [n]} m_{v^*_t}(A) \, \big(1-F(\mathbf{e}^{t,b}_{[n]\setminus A})\big) = 1-\sum_{A\subseteq [n]} m_{v^*_t}(A) \,
F(\mathbf{e}^{t,b}_{[n]\setminus A}).\label{eq:RSmvsWLP}
\end{eqnarray}
\end{theorem}

As far as the mean time-to-failure of the system is concerned, from (\ref{eq:RSmvWLP}) and (\ref{eq:RSmvsWLP}) we immediately obtain
\begin{eqnarray}
\mathrm{MTTF}_S &=& \sum_{A\subseteq [n]}\sum_{B\subseteq A} (-1)^{|A|-|B|}\, \int_0^{w(B)} R(\mathbf{e}^{0,t}_A)\, dt,\label{eq:MTTFsWB}\\
\mathrm{MTTF}_S &=& \sum_{A\subseteq [n]}\sum_{B\subseteq A} (-1)^{|A|-|B|}\, \int_{w([n]\setminus B)}^{\infty}
\big(1-F(\mathbf{e}^{t,\infty}_{[n]\setminus A})\big)\, dt.\label{eq:MTTFsWB2}
\end{eqnarray}

When the variables $T_1,\ldots,T_n$ are independent, from (\ref{eq:RSvWLP}) we immediately retrieve the formula (see Marichal \cite{Mar08}):
$$
R_S(t) = \sum_{A\subseteq [n]} v_t(A) \prod_{i\in A}R_i(t)\prod_{i\in [n]\setminus A}(1-R_i(t)).
$$
Considering the family $\{\overline{\phi}_{v_t}\, : \, t\geqslant 0\}$, where $\overline{\phi}_{v_t}$ is the multilinear extension of
$\phi_{v_t}$, we then observe that
\begin{equation}\label{eq:R-Ind-MLE-WLP}
R_S(t)=\overline{\phi}_{v_t}(R_1(t),\ldots,R_n(t))
\end{equation}
and $\overline{\phi}_{v_t}$ can be chosen from among the forms given in Table~\ref{tab:StrFun}, where each $v$ should be replaced with $v_t$.
Also, from (\ref{eq:MTTFsWB}) and (\ref{eq:MTTFsWB2}) we immediately derive
\begin{eqnarray}
\mathrm{MTTF}_S &=& \sum_{A\subseteq [n]}\sum_{B\subseteq A} (-1)^{|A|-|B|}\, \int_0^{w(B)} \prod_{i\in A} R_i(t)\, dt,\label{eq:MTTFsWB-Ind}\\
\mathrm{MTTF}_S &=& \sum_{A\subseteq [n]}\sum_{B\subseteq A} (-1)^{|A|-|B|}\, \int_{w([n]\setminus B)}^{\infty} \coprod_{i\in A} R_i(t)\,
dt.\nonumber
\end{eqnarray}

Let us now examine some examples by considering special w.l.p.\ functions. They generalize the classical examples considered in \S\ref{sec:ex5}
(series, parallel, and $k$-out-of-$n$ structures).
\begin{enumerate}
\item \emph{Weighted minimum.} A weighted minimum function is a w.l.p.\ function $p_w:L^n\to L$ whose underlying set function $w:2^{[n]}\to L$
fulfills
$$
w([n]\setminus (A\cup B))=w([n]\setminus A)\wedge w([n]\setminus B).
$$
Such a function fulfills equation (\ref{eq:Minitivity}) and is of the form (see \cite[\S{5.2}]{Mar00c})
$$
p_w(\mathbf{t})=\bigwedge_{i=1}^n\big(w(\{i\})\vee t_i\big).
$$
It then corresponds to a series structure with a lower bound on the lifetime of each component. By using (\ref{eq:RSmvsWLP}), we can easily show
that
$$
R_S(t)=1-\sum_{\textstyle{A\subseteq [n]\atop A\neq\varnothing}}(-1)^{|A|+1}F(\mathbf{e}_{[n]\setminus A}^{t,\infty})\prod_{i\in
A}(1-v_t(\{i\}))
$$
and, in case of independence (see~(\ref{eq:R-Ind-MLE-WLP})),
$$
R_S(t)=\prod_{i=1}^n \big(v_t(\{i\})\amalg R_i(t)\big).
$$

\item \emph{Weighted maximum.} A weighted maximum function is a w.l.p.\ function $p_w:L^n\to L$ whose underlying set function $w:2^{[n]}\to L$
fulfills
$$
w(A\cup B)=w(A)\vee w(B).
$$
Such a function fulfills equation (\ref{eq:Maxitivity}) and is of the form (see \cite[\S{5.2}]{Mar00c})
$$
p_w(\mathbf{t})=\bigvee_{i=1}^n\big(w(\{i\})\wedge t_i\big).
$$
It then corresponds to a parallel structure with an upper bound on the lifetime of each component. By using (\ref{eq:RSmvWLP}), it is also
straightforward to show that
$$
R_S(t)=\sum_{\textstyle{A\subseteq [n]\atop A\neq\varnothing}}(-1)^{|A|+1}R(\mathbf{e}_{A}^{0,t})\prod_{i\in A}v_t(\{i\})
$$
and, in case of independence (see~(\ref{eq:R-Ind-MLE-WLP})),
$$
R_S(t)=\coprod_{i=1}^n v_t(\{i\})\, R_i(t).
$$

\item \emph{Symmetric w.l.p.\ function.} We can generalize the $k$-out-of-$n$ type structures simply by considering symmetric w.l.p.\ functions
$p_w:L^n\to L$. The underlying set functions are cardinality based, i.e., such that $w(A)=w(B)$ whenever $|A|=|B|$. If we define the function
$\widetilde{w}:\{0,1,\ldots,n\}\to L$ by $w(A)=\widetilde{w}(|A|)$, we can easily show \cite{DukMar08} that any symmetric w.l.p.\ function can
always be put in the form
$$
p_w(\mathbf{t})=\bigvee_{k=0}^n\big(\widetilde{w}(k)\wedge f_{n-k+1}(\mathbf{t})\big),
$$
where $f_{n-k+1}$ is the order statistic function defining the $k$-out-of-$n$ structure (see \S\ref{sec:ex5}). Moreover, we can show that
$$
R_S(t)=\Pr(f_{n-k(t)+1}(\mathbf{T})>t)=\Pr(|\mathbf{X}(t)|\geqslant k(t)),
$$
where $k(t):=\min\{k,n+1 : \widetilde{w}(k)>t\}$, which generalizes (\ref{eq:SymPP}). For more details, see Dukhovny and
Marichal~\cite{DukMar08}.
\end{enumerate}

\begin{example}\label{ex:EXPLE-WLP}
Let us calculate $\mathrm{MTTF}_S$ when $R_i(t)=e^{-\lambda_i t}$ $(i=1,\ldots,n)$. Using (\ref{eq:MTTFsWB-Ind}) and setting
$\lambda_A=\sum_{i\in A}\lambda_i$, we obtain (see \cite{Mar08})
\begin{eqnarray*}
\mathrm{MTTF}_S &=& \sum_{A\subseteq [n]} \sum_{B\subseteq A} (-1)^{|A|-|B|} \int_0^{w(B)} e^{-\lambda_A t}\, dt\\
&=& w(\varnothing)+\sum_{\textstyle{A\subseteq [n]\atop A\neq\varnothing}} \sum_{B\subseteq A} (-1)^{|A|-|B|} \frac{1-e^{-\lambda_A
w(B)}}{\lambda_A}\, .
\end{eqnarray*}
When $p_w$ is an l.p.\ function (that is, $w(A)\in\{0,\infty\}$ and $w(\varnothing)=0$), we retrieve formula (\ref{eq:MTTFSExp}). Indeed,
$$
\mathrm{MTTF}_S = \sum_{\textstyle{A\subseteq [n]\atop A\neq\varnothing}} \frac{1}{\lambda_A} \sum_{\textstyle{B\subseteq A\atop w(B)=\infty}}
(-1)^{|A|-|B|} = \sum_{\textstyle{A\subseteq [n]\atop A\neq\varnothing}} \frac{m_v(A)}{\lambda_A}\, .
$$
\end{example}

\section{Distribution functions of w.l.p.\ functions}

The articles \cite{Duk07,Mar06,Mar08} on which this paper is partly based were motivated by the exact computation of the distribution functions
and the moments of l.p.\ functions and w.l.p.\ functions of random variables.

In this final section we point out the fact that calculating the distribution function of an arbitrary w.l.p.\ function amounts to calculating
the reliability function of a semicoherent system with possible lower and upper bounds on component lifetimes.

Let $L\subseteq \overline{\R}$ be a totally ordered lattice, let $p_w:L^n\to L$ be a w.l.p., and let $T_1,\ldots,T_n$ be $L$-valued random
variables.

The distribution function of the random variable $p_w(\mathbf{T})$ is defined as
$$
F_{p_w}(t)=\Pr(p_w(\mathbf{T})\leqslant t) \qquad (t\in L).
$$
Clearly, this function fulfills the identity $F_{p_w}(t)=1-R_S(t)$, where $R_S(t)$ is the reliability function of the coherent system described
by the w.l.p.\ function $p_w$.

Using formulas (\ref{eq:RSvWLP})--(\ref{eq:RSmvsWLP}), we then obtain immediately the following formulas for $F_{p_w}(t)$:
\begin{eqnarray*}
F_{p_w}(t) &=& 1-\sum_{A\subseteq [n]} v_t(A)\,\Pr(\mathbf{X}(t)=\mathbf{e}^{0,1}_A)\\
F_{p_w}(t) &=& \sum_{A\subseteq [n]} v_t^*(A)\,\Pr(\mathbf{X}(t)=\mathbf{e}^{0,1}_{[n]\setminus A})\\
F_{p_w}(t) &=& 1-\sum_{A\subseteq [n]} m_{v_t}(A) \, R(\mathbf{e}^{a,t}_A)\\
F_{p_w}(t) &=& \sum_{A\subseteq [n]} m_{v^*_t}(A) \, F(\mathbf{e}^{t,b}_{[n]\setminus A}),
\end{eqnarray*}
where $v_t(A)=\Ind(w(A)>t)$ and $v^*_t(A)=1-\Ind(w([n]\setminus A)>t)$.

When the arguments $T_1,\ldots,T_n$ are independent, each $T_i$ having distribution function $F_i(t)$, we obtain (see Marichal~\cite{Mar08}):
\begin{eqnarray*}
F_{p_w}(t) &=& 1-\sum_{A\subseteq [n]} v_t(A)\,\prod_{i\in A}(1-F_i(t))\prod_{i\in [n]\setminus A} F_i(t)\\
F_{p_w}(t) &=& \sum_{A\subseteq [n]} v_t^*(A)\,\prod_{i\in A}F_i(t)\prod_{i\in [n]\setminus A} (1-F_i(t))\\
F_{p_w}(t) &=& 1-\sum_{A\subseteq [n]} m_{v_t}(A) \, \prod_{i\in A}(1-F_i(t))\\
F_{p_w}(t) &=& \sum_{A\subseteq [n]} m_{v^*_t}(A) \, \prod_{i\in A}F_i(t).
\end{eqnarray*}

\section{Conclusion}

We have discussed the formal parallelism between two representations of semicoherent systems: structure functions and l.p.\ functions. Their
languages are shown to be equivalent in many ways. The l.p.\ language is demonstrated to have significant advantages. One is the natural
generalization to w.l.p.\ functions and corresponding systems with bounded subsystem lifetimes. The other is
the fact that, due to the distributive property of the indicator function $\mathrm{Ind}(\cdot)$ with respect to lattice operations (see proofs
of Theorems~\ref{thm:SFvsLP} and \ref{thm:SFvsWLP}), the l.p.\ description is a very natural tool to connect the system's structure to the
lattice of typical reliability events of the kind $T\leqslant t$, to connect the system's purpose, as encoded in the l.p.\ function, to the
system's equipment, as expressed in the joint distribution of units' lifetimes.



\end{document}